\documentclass[preprint,review,12pt]{elsarticle}
\usepackage{amsfonts}
\usepackage[fleqn]{amsmath}
\usepackage{epsfig}
\usepackage{subfigure}
\usepackage{amsmath}
\usepackage{tikz}
\usepackage[latin1]{inputenc}
\usepackage{verbatim}
\usepackage{mathptmx}
\usepackage{graphicx}
\usepackage{dcolumn}
\usepackage{bm}

\usepackage{mathrsfs}
\usepackage[colorlinks,citecolor=blue,urlcolor=blue]{hyperref}

\numberwithin{equation}{section}
\numberwithin{figure}{section}
\numberwithin{table}{section}
\newtheorem{thm}{Theorem}[section]

\newtheorem{definition}{Definition}[section]
\newtheorem{corr}{Corollary}[section]
\newdefinition{rmk}{Remark}[section]
\newproof{pr}{Proof}
\newdefinition{scheme}{Scheme}

\usepackage{lineno}

\journal{XX}

\begin{document}

\begin{frontmatter}

\title{Stochastic multi-symplectic Runge-Kutta methods for stochastic Hamiltonian PDEs}

\author[rvt]{Liying Zhang}
\ead{lyzhang@lsec.cc.ac.cn}
\author[els]{Lihai Ji\corref{cor1}}
\ead{jilihai@lsec.cc.ac.cn}

\address[rvt]{School of Mathematical Science, China University of Mining and Technology, Beijing 100083, China}
\address[els]{Institute of Applied Physics and Computational Mathematics, Beijing 100094, China}

\cortext[cor1]{Corresponding author}
\fntext[fn1]{The first author is supported by NNSFC (NO. 11601514 and NO. 11771444) and the second author is supported by NNSFC (NO. 11471310 and NO. 11601032).}

\begin{abstract}
In this paper, we consider stochastic Runge-Kutta methods
for stochastic Hamiltonian partial differential equations and present some
sufficient conditions for multi-symplecticity of stochastic Runge-Kutta methods of stochastic Hamiltonian
partial differential equations. Particularly, we apply these ideas to stochastic Maxwell equations with multiplicative noise, possessing the stochastic multi-symplectic conservation law and energy conservation law. Theoretical analysis shows that the methods can preserve both the discrete stochastic multi-symplectic conservation law and discrete energy conservation law almost surely.
\end{abstract}

\end{frontmatter}


\section{Introduction}\label{1.1}
Consider the following stochastic partial differential equations (SPDEs) in the sense of Stratonovich
\begin{equation}\label{SHPDE}
{\rm {\bf K}}{\rm d}_tz+{\rm {\bf L}}z_{x}{\rm d}t=\nabla_z S_1(z){\rm d}t+\nabla_z S_2(z)\circ{\rm d}W(t),~~z\in\mathbb{R}^n,~n\geq2,
\end{equation}
where ${\rm {\bf K}}$, ${\rm {\bf L}}$ are skew-symmetric matrices, and $S_1$, $S_2$ are
smooth functions of state variable $z$. Let $H=L^2(\mathbb{R},\mathbb{R})$ and
let $(\Omega,\mathcal{F},P)$ be a probability space with a normal filtration $\{\mathcal{F}_{t}\}_{0\leq t\leq T}$.
Moreover, let $W:[0,T]\times\Omega\rightarrow H$ be a standard $Q$-Wiener process with
respect to $\{{\mathcal F}_{t}\}_{0\leq t\leq T}$, with a trace class operator $Q: H\rightarrow H$.
The SPDEs are called the stochastic Hamiltonian PDEs (see \cite{JWH2013}). Moreover, the stochastic Hamiltonian PDEs \eqref{SHPDE} possess the stochastic multi-symplectic conservation law (Theorem 2.2 in \cite{JWH2013}) as follows
\begin{equation}\label{s_multi_symplectic_law}
  {\rm d}_t\omega(t,x)+\partial_x\kappa(t,x){\rm d}t=0,~a.s.,
\end{equation}
where $\omega(t,x)=\frac{1}{2}{\rm d}z\wedge {\rm {\bf K}}{\rm d}z$, $\kappa(t,x)=\frac{1}{2}{\rm d}z\wedge {\rm {\bf L}}{\rm d}z$ are differential 2-forms associated with the two skew-symmetric matrices ${\rm {\bf K}}$ and ${\rm {\bf L}}$, respectively.

Many physical and engineering phenomena are modeled by stochastic Hamiltonian PDEs. The advantage of modeling using these so-called stochastic Hamiltonian PDEs is that stochastic Hamiltonian PDEs possess stochastic multi-symplectic geometric structure \eqref{s_multi_symplectic_law} and then are able to more fully capture the behavior of interesting phenomena. Therefore, it requires the corresponding numerical schemes could exactly preserve the discrete stochastic multi-symplectic structure, which are known as stochastic multi-symplectic schemes. In recent years, many researchers have studied different stochastic Hamiltonian PDEs and various stochastic multi-symplectic numerical schemes have been developed, analyzed and tested, see, stochastic Maxwell equaitons \cite{CHZ2016,HJZ2014,HJZC2017}; stochastic nonlinear Schr\"{o}dinger equation \cite{CHJK2017,CLHZ2017,HWZ2017,JWH2013}; stochastic Korteweg-de Vries equation \cite{JWH2012}.

In the well developed numerical analysis of deterministic ODEs, derivative free approximation methods are of particular interest. Especially Runge-Kutta (RK) methods are widely used. In the last decade, it has been widely recognized that the stochastic RK methods play an important role in numerically solving stochastic ODEs and popularly employed, see, e.g.,~\cite{BB2001,BB2014,HXW2015,KP1999,MD2012,MD2015,Rosser2009,TV2002,Wang2008} and the references therein. The symplectic condition of stochastic RK methods for stochastic Hamiltonian ODEs was obtained firstly in \cite{MD2012}. Furthermore, \cite{MD2015} derived the symplectic conditions of stochastic partitioned Runge-Kutta (PRK) methods. Recently, based on variational principle, \cite{CH2016} proposed a general class of stochastic symplectic RK methods in the temporal direction to the stochastic Schr\"{o}dinger equation in the Stratonovich sense and showed that the methods preserve the charge conservation law. In particular, the authors present the mean-square convergence order of the semidiscrete scheme is 1 under
appropriate assumptions.

For deterministic Hamiltonian PDEs, the multi-symplecticity of RK methods has been studied in \cite{HLS2005}. To the best of our knowledge, there has been no work in the literature which studies the multi-symplectic RK methods for stochastic Hamiltonian PDEs \eqref{SHPDE}. Motivated by \cite{CH2016,MD2012,MD2015}, we consider the general case of stochastic Hamiltonian PDEs, investigate the multi-symplecticity of stochastic RK methods, and present some sufficient conditions for stochastic multi-symplectic RK methods in this paper. To this end, we main utilize the corresponding variation form of stochastic Hamiltonian PDEs and the chain rule of Stratonovich integral, then derive the multi-symplectic conditions by a tedious computation.

The rest of this paper is organized as follows. In Section 2, we apply stochastic RK methods to solve \eqref{SHPDE} and then present conditions for multi-symplecticity of stochastic RK methods. In particular, an one-stage stochastic multi-symplectic RK method is given. In Section 3, the multi-symplecticity of stochastic RK methods for the three-dimensional stochastic Maxwell equations with multiplicative noise is discussed. Concluding remarks are presented in Section 4.

\section{Stochastic multi-symplectic RK methods}
\subsection{Stochastic multi-symplectic conditions}
In this section, we give the conditions of multi-symplecticity
of stochastic RK discretization for general stochastic Hamiltonian PDEs. In the sequel, we denote $Z_{m}^{k}\approx z(c_m h,d_k\tau)$, $Z_{m}^{1}\approx z(c_m h,\tau)$, $Z_{1}^{k}\approx z(h,d_k\tau)$, $c_m=\sum_{n=1}^s\widetilde{a}_{mn}$, $d_k=\sum_{j=1}^ra_{kj}$. To solve \eqref{SHPDE}, we present a class of stochastic RK methods with r-stage in temporal direction and s-stage in spatial direction, respectively, namely
\begin{subequations}\label{PRK}
\begin{align}
  &Z_{m}^{k}=z_{m}^{p}+\tau\sum_{j=1}^{r}a_{kj}\delta_t^{m,j}Z_{m}^{j},\quad \forall~p=0,1,\cdots,P\label{PRK_sub1}\\[1mm]
  &z_{m}^{p+1}=z_{m}^{p}+\tau\sum_{k=1}^{r}b_{k}\delta_t^{m,k}Z_{m}^{k},\quad \forall~p=0,1,\cdots,P\label{PRK_sub2}\\[1mm]
  &Z_{m}^{k}=z_{i}^{k}+h\sum_{n=1}^{s}\widetilde{a}_{mn}\delta_{x}^{n,k}Z_{n}^{k},\quad \forall~i=0,1,\cdots,I\label{PRK_sub3}\\[1mm]
  &z_{i+1}^{k}=z_{i}^{k}+h\sum_{m=1}^{s}\widetilde{b}_{m}\delta_{x}^{m,k}Z_{m}^{k},\quad \forall~i=0,1,\cdots,I\label{PRK_sub4}\\[1mm]
  &\tau{\rm {\bf K}}\delta_{t}^{m,k}Z_{m}^{k}+\tau{\rm {\bf L}}\delta_{x}^{m,k}Z_{m}^{k}=\tau\nabla_{z}S_1(Z_{m}^{k})+\nabla_{z}S_2(Z_{m}^{k})\Delta W_m^k,\label{PRK_sub5}
  \end{align}
\end{subequations}
where $\tau$, $h$ are the stepsizes in time and spatial directions, respectively. $\delta_t^{m,k}$ and $\delta_{x}^{m,k}$ the discretizations of the partial derivatives $\partial_t$ and $\partial_x$, respectively.  The increment $\Delta W_m^k:=W(t_{k+1},x_m,\omega)-W(t_{k},x_m,\omega)$ is given by
\begin{equation}\label{Q_process}
  \Delta W_m^k(\omega):=\sum_{j=1}^{\infty}\sqrt{\eta_j}e_j(x_m)\left(\beta_j(t^{k+1},\omega)-\beta_j(t^k,\omega)\right)
\end{equation}
for all $\omega\in\Omega$. Here $e_j,~j\in\mathbb{N}$ is an orthonormal basis of $H$ consisting of eigenfunctions of $Q$ such that $Qe_j=\eta_j e_j,~j\in\mathbb{N}$ and $\beta_j,~j\in\mathbb{N}$ are independent real-valued Brownian motions on the
probability space $(\Omega,\mathcal{F},P,\{\mathcal{F}_{t}\}_{0\leq t\leq T})$.

In order not to complicate the notation, and for sake of brevity, we shall use the abbreviations, $\delta_{t}$ denotes $\delta_{t}^{m,k}$, $\delta_{x}$ denotes $\delta_{x}^{m,k}$ and so on.
\begin{definition}
  A numerical method with the approximating solution $z_m^{k}$ for  \eqref{SHPDE} is said to be stochastic multi-symplectic if it satisfies
\begin{equation}\label{def_sm}
  \tau\delta_{t}\omega_{m}^{k}+h\delta_{x}\kappa_{m}^{k}=0,~a.s.,
\end{equation}
where $\omega_m^k=\frac{1}{2}{\rm d}z_{m}^{k}\wedge {\rm {\bf K}}{\rm d}z_{m}^{k}$, $\kappa_{m}^{k}=\frac{1}{2}{\rm d}z_{m}^{k}\wedge {\rm {\bf L}}{\rm d}z_{m}^{k}$.
\end{definition}


As follows we state the main theorem of this paper.
\begin{thm}
Assume that the coefficients $a_{kj},\widetilde{a}_{mn},b_k,\widetilde{b}_m$ of \eqref{PRK} satisfy the relations
\begin{align}\label{cond0}
  b_kb_j-b_ka_{kj}-b_ja_{jk}=0~~{\rm and}~~\widetilde{b}_m\widetilde{b}_n-\widetilde{b}_m\widetilde{a}_{mn}-\widetilde{b}_n\widetilde{a}_{nm}=0,
\end{align}
for all $k,j=1,\cdots,r$ and $n,m=1,\cdots,s$, then the stochastic RK methods \eqref{PRK} is stochastic multi-symplectic with the discrete stochastic multi-symplectic conservation law almost surely
\begin{equation}\label{def_sm_prk}
  \frac{\omega^{p+1}-\omega^{p}}{\tau}+\frac{\kappa_{i+1}-\kappa_{i}}{h}=0,\quad \forall~p=0,1,\cdots,P;~i=0,1,\cdots,I,
\end{equation}
where
\begin{equation*}
\begin{split}
\omega^p&=\frac{1}{2}\sum_{m=1}^s\widetilde{b}_m{\rm d}z_{m}^{p}\wedge {\rm {\bf K}}{\rm d}z_{m}^{p},\quad
\kappa_{i}=\frac{1}{2}\sum_{k=1}^rb_k{\rm d}z_{i}^{k}\wedge {\rm {\bf L}}{\rm d}z_{i}^{k}.
\end{split}
\end{equation*}
\end{thm}

\begin{pr}
Differentiating \eqref{PRK_sub1} and \eqref{PRK_sub2}, it holds
\begin{subequations}
  \begin{align}
  {\rm d}Z_{m}^{k}&={\rm d}z_{m}^{p}+\tau\sum_{j=1}^{r}a_{kj}{\rm d}(\delta_tZ_{m}^{j}),\label{Diff_1}\\[1mm]
  {\rm d}z_{m}^{p+1}&={\rm d}z_{m}^{p}+\tau\sum_{k=1}^{r}b_{k}{\rm d}(\delta_tZ_{m}^{k}).\label{Diff_2}
  \end{align}
\end{subequations}
Then we have
\begin{equation}\label{Proof_1}
  \begin{split}
    {\rm d}z_m^{p+1}\wedge{\rm {\bf K}}{\rm d}z_m^{p+1}&=\left({\rm d}z_{m}^{p}+\tau\sum_{k=1}^{r}b_{k}{\rm d}\left(\delta_tZ_{m}^{k}\right)\right)\wedge{\rm {\bf K}}\left({\rm d}z_{m}^{p}+\tau\sum_{k=1}^{r}b_{k}{\rm d}\left(\delta_tZ_{m}^{k}\right)\right)\\
    &={\rm d}z_{m}^{p}\wedge{\rm {\bf K}}{\rm d}z_{m}^{p}+\tau\sum_{k=1}^rb_k\left({\rm d}z_{m}^{p}\wedge{\rm {\bf K}}{\rm d}(\delta_tZ_{m}^{k})+{\rm d}(\delta_tZ_{m}^{k})\wedge{\rm {\bf K}}{\rm d}z_{m}^{p}\right)\\
    &\quad\quad+\tau^2\sum_{k=1}^r\sum_{j=1}^rb_kb_j\left({\rm d}(\delta_tZ_{m}^{k})\wedge{\rm {\bf K}}{\rm d}(\delta_tZ_{m}^{j})\right).
  \end{split}
\end{equation}

Substituting \eqref{Diff_1} into the second term of the right-side of the equation \eqref{Proof_1}, it yields
\begin{equation*}
  \begin{split}
    {\rm d}z_m^{p+1}\wedge{\rm {\bf K}}{\rm d}z_m^{p+1}&={\rm d}z_{m}^{p}\wedge{\rm {\bf K}}{\rm d}z_{m}^{p}+\tau\sum_{k=1}^rb_k\left({\rm d}Z_m^k\wedge{\rm {\bf K}}{\rm d}(\delta_tZ_{m}^{k})+{\rm d}(\delta_tZ_{m}^{k})\wedge{\rm {\bf K}}{\rm d}Z_m^k\right)\\
    &\quad\quad+\tau^2\sum_{k=1}^r\sum_{j=1}^r(b_kb_j-b_ka_{kj}-b_ja_{jk})\left({\rm d}(\delta_tZ_{m}^{k})\wedge{\rm {\bf K}}{\rm d}(\delta_tZ_{m}^{j})\right).
  \end{split}
\end{equation*}
Using \eqref{cond0} and the skew-symmetry of matrix ${\rm {\bf K}}$, we have
\begin{equation}\label{Proof_2}
  \begin{split}
    {\rm d}z_m^{p+1}\wedge{\rm {\bf K}}{\rm d}z_m^{p+1}&={\rm d}z_{m}^{p}\wedge{\rm {\bf K}}{\rm d}z_{m}^{p}\\
    &\quad\quad+\tau\sum_{k=1}^rb_k\left({\rm d}Z_m^k\wedge{\rm {\bf K}}{\rm d}(\delta_tZ_{m}^{k})+{\rm d}(\delta_tZ_{m}^{k})\wedge{\rm {\bf K}}{\rm d}Z_m^k\right)\\
    &={\rm d}z_{m}^{p}\wedge{\rm {\bf K}}{\rm d}z_{m}^{p}+2\tau\sum_{k=1}^rb_k{\rm d}Z_m^k\wedge{\rm {\bf K}}{\rm d}(\delta_tZ_{m}^{k}).
  \end{split}
\end{equation}

Similarly, it can obtain
\begin{equation}\label{Proof_3}
  \begin{split}
    {\rm d}z_{i+1}^{k}\wedge{\rm {\bf L}}{\rm d}z_{i+1}^{k}&={\rm d}z_{i}^{k}\wedge{\rm {\bf L}}{\rm d}z_{i}^{k}\\
    &\quad\quad+h\sum_{m=1}^s\widetilde{b}_m\left({\rm d}Z_m^k\wedge{\rm {\bf L}}{\rm d}(\delta_xZ_{m}^{k})+{\rm d}(\delta_xZ_{m}^{k})\wedge{\rm {\bf L}}{\rm d}Z_m^k\right)\\
    &={\rm d}z_{i}^{k}\wedge{\rm {\bf L}}{\rm d}z_{i}^{k}+2h\sum_{m=1}^s\widetilde{b}_m{\rm d}Z_m^k\wedge{\rm {\bf L}}{\rm d}(\delta_xZ_{m}^{k}).
  \end{split}
\end{equation}
Multiplying both sides of \eqref{Proof_2} and \eqref{Proof_3} with $\frac{1}{2}h\widetilde{b}_m$, $\frac{1}{2}\tau b_k$, summing over all spatial grid
points $m$ and temporal grid points $k$, respectively, and adding them together, we haves
\begin{equation}\label{Proof_4}
  \begin{split}
    \frac{h}{2}\sum_{m=1}^{s}&\widetilde{b}_m{\rm d}z_{m}^{p+1}\wedge{\rm {\bf K}}{\rm d}z_{m}^{p+1}+\frac{\tau}{2}\sum_{k=1}^{r}b_k{\rm d}z_{i+1}^{k}\wedge{\rm {\bf L}}{\rm d}z_{i+1}^{k}\\
    &=\frac{h}{2}\sum_{m=1}^{s}\widetilde{b}_m{\rm d}z_{m}^{p}\wedge{\rm {\bf K}}{\rm d}z_{m}^{p}+\frac{\tau}{2}\sum_{k=1}^{r}b_k{\rm d}z_{i}^{k}\wedge{\rm {\bf L}}{\rm d}z_{i}^{k}\\
    &+\tau h\sum_{m=1}^s\sum_{k=1}^{r}b_k\widetilde{b}_m\left({\rm d}Z_m^k\wedge{\rm {\bf K}}{\rm d}(\delta_tZ_{m}^{k})+{\rm d}Z_m^k\wedge{\rm {\bf L}}{\rm d}(\delta_xZ_{m}^{k})\right).
  \end{split}
\end{equation}
Noticing that, if we take the differential in the phase space on both sides of \eqref{PRK_sub5}, we can deduce that
\begin{equation*}
  {\rm {\bf K}}{\rm d}(\delta_{t}Z_{m}^{k})+{\rm {\bf L}}{\rm d}(\delta_{x}Z_{m}^{k})=D_{zz}S_1(Z_{m}^{k}){\rm d}Z_m^k+D_{zz}S_2(Z_{m}^{k}){\rm d}Z_m^k\frac{\Delta W_m^k}{\tau}
\end{equation*}
then we use ${\rm d}Z_m^k$ to perform wedge product with the above equation, it yields
\begin{equation}\label{Proof_5}
  {\rm d}Z_m^k\wedge{\rm {\bf K}}{\rm d}(\delta_{t}Z_{m}^{k})+{\rm d}Z_m^k\wedge{\rm {\bf L}}{\rm d}(\delta_{x}Z_{m}^{k})=0,
\end{equation}
where the equality is due to the symmetry of $D_{zz}S_1(Z_{m}^{k})$ and $D_{zz}S_2(Z_{m}^{k})$. Combining \eqref{Proof_4} and \eqref{Proof_5}, we get the stochastic multi-symplectic conservation law \eqref{def_sm_prk}.
\end{pr}

\subsection{One-stage stochastic multi-symplectic RK method}
In this subsection, we use the stochastic multi-symplectic conditions \eqref{cond0} to construct an one-stage stochastic multi-symplectic RK method. Consider one-stage stochastic multi-symplectic RK methods in the following form
\begin{center}
\begin{tabular}{c|c}
       & a \\
\hline & b
\end{tabular},
~~~~\begin{tabular}{c|c}
       & $\widetilde{a}$ \\
\hline & $\widetilde{b}$
\end{tabular}.
\end{center}
Using the stochastic multi-symplectic conditions \eqref{cond0}, the following results hold:
\begin{align}
  b=2a,~ \widetilde{b}=2\widetilde{a}.
\end{align}
In particular, choosing $b=\widetilde{b}=1$, we get the scheme as
\begin{center}
\begin{tabular}{c|c}
       & {\rm 1/2} \\
\hline & {\rm 1}
\end{tabular},
~~~~\begin{tabular}{c|c}
       & {\rm 1/2} \\
\hline & {\rm 1}
\end{tabular},
\end{center}
more precisely,
\begin{equation}
  {\rm {\bf K}}\left(\frac{z_{i+\frac{1}{2}}^{p+1}-z_{i+\frac{1}{2}}^{p}}{\tau}\right)+{\rm {\bf L}}\left(\frac{z_{i+1}^{p+\frac{1}{2}}-z_{i}^{p+\frac{1}{2}}}{h}\right)=\nabla_{z}S_1(z_{i+\frac{1}{2}}^{p+\frac{1}{2}})+\nabla_{z}S_2(z_{i+\frac{1}{2}}^{p+\frac{1}{2}})\frac{\Delta W_{i+\frac{1}{2}}^{p+\frac{1}{2}}}{\tau},
\end{equation}
which is equivalent to implicit midpoint scheme in \cite{JWH2013}.

\section{Application of stochastic multi-symplectic RK methods}
In this section, we apply the stochastic RK methods \eqref{PRK} to the stochastic Maxwell equations with multiplicative noise in statistical radiophysics \cite{RKT1987}. We obtain the sufficient conditions for multi-symplecticity of stochastic RK methods. Moreover, we derive the discrete energy conservation law under stochastic multi-symplectic RK methods.

Considering the following stochastic Maxwell equations with multiplicative noise in Stratonovich sense \cite{HJZC2017,LSY2010,RKT1987}
\begin{equation}\label{stochastic maxwell equations}
\begin{split}
{\rm d}\mathcal{E}(x,t)&=A_M\mathcal{E}(x,t){\rm d}t+B(\mathcal{E})(x,t)\circ {\rm d}W(t),~x\in D,~t>0\\
\mathcal{E}_0(x)&=\xi,~x\in D
\end{split}
\end{equation}
in $L^{2}(D)^6=L^{2}(D)^3\times L^{2}(D)^3,~D\subset \mathbb{R}^3$ driven by a standard $Q$-Wiener process $W(t)$ and the perfect conductor boundary condition $n\times E=0$ on $\partial D$. Let $\mathcal{E}=\left(H^T,E^T\right)^T$, $H=(H_1,H_2,H_3)$ and $E=(E_1,E_2,E_3)$, then the Maxwell operator is given by
\begin{equation*}
  A_M\begin{pmatrix}
H\\
E
\end{pmatrix}=\begin{pmatrix}
0&-\nabla\times\\
\nabla\times&0
\end{pmatrix}\begin{pmatrix}
H\\
E
\end{pmatrix},
\end{equation*}
and
\begin{equation*}
  B(\mathcal{E})=\begin{pmatrix}
0&\lambda\\
-\lambda &0
\end{pmatrix}\begin{pmatrix}
H\\
E
\end{pmatrix}.
\end{equation*}

In \cite{HJZC2017}, the authors show that stochastic Maxwell equations can be written
\begin{equation}\label{12}
{\rm {\bf K}}{\rm d}_{t}\mathcal{E}+{\rm {\bf L}}_{1}\mathcal{E}_{x}{\rm d}t+{\rm {\bf L}}_{2}\mathcal{E}_{y}{\rm d}t+{\rm {\bf L}}_{3}\mathcal{E}_{z}{\rm d}t=\nabla S(\mathcal{E})\circ
{\rm d}W(t),
\end{equation}
with
\begin{equation}\label{SS}
S(\mathcal{E})=\frac{\lambda}{2}\left(|E_1|^{2}+|E_2|^{2}+|E_3|^{2}+|H_1|^{2}+|H_2|^{2}+|H_3|^{2}\right)
\end{equation}
and
\begin{equation*}
{\rm {\bf K}}=\left(\begin{array}{ccccccc}
0&-I_{3\times3}\\[1mm]
I_{3\times3}&0\\[1mm]
\end{array}\right),\quad
{\rm {\bf L}}_{i}=\left(\begin{array}{cccccccc}
\mathcal{D}_{i}&0\\[1mm]
0&\mathcal{D}_{i}\\[1mm]
\end{array}
\right),~~i=1, 2, 3.
\end{equation*}
The sub-matrix $I_{3\times3}$ is a $3\times3$ identity matrix and
\begin{small}
	$$
	\mathcal{D}_{1}=\left(\begin{array}{ccccccc}
	0&0&0\\[1mm]
	0&0&-1\\[1mm]
	0&~1&0\\[1mm]
	\end{array}\right),
	\mathcal{D}_{2}=\left(\begin{array}{ccccccc}
	0&0&~1\\[1mm]
	0&0&0\\[1mm]
	-1&0&0\\[1mm]
	\end{array}\right),
	\mathcal{D}_{3}=\left(\begin{array}{ccccccc}
	0&-1&0\\[1mm]
	~1&0&0\\[1mm]
	0&0&0\\[1mm]
	\end{array}\right).
	$$
	\end{small}

Now we investigate the stochastic multi-symplecticity of RK methods for the equations \eqref{stochastic maxwell equations}. The stochastic RK method \eqref{PRK} applied to the equation \eqref{stochastic maxwell equations} is
\begin{subequations}\label{RK00}
\begin{align}
  &\Upsilon_{mplk}=\mathcal{E}_{mpl}^{\rho}+\tau\sum_{j=1}^{r}a_{kj}\delta_t\Upsilon_{mplj},\quad \forall~\rho=0,1,\cdots,N\label{GME_sub1}\\[1mm]
  &\mathcal{E}_{mpl}^{\rho+1}=\mathcal{E}_{mpl}^{\rho}+\tau\sum_{k=1}^{r}b_{k}\delta_t\Upsilon_{mplk},\quad \forall~\rho=0,1,\cdots,N\label{GME_sub2}\\[1mm]
  &\Upsilon_{mplk}=\mathcal{E}_{i_1pl}^{k}+\Delta x\sum_{n=1}^{s}\widetilde{a}_{mn}\delta_{x}\Upsilon_{nplk},\quad \forall~i_1=0,1,\cdots,I_1\label{GME_sub3}\\[1mm]
  &\mathcal{E}_{(i_1+1)pl}^{k}=\mathcal{E}_{i_1pl}^{k}+\Delta x\sum_{m=1}^{s}\widetilde{b}_{m}\delta_{x}\Upsilon_{mplk},\quad \forall~i_1=0,1,\cdots,I_1\label{GME_sub4}\\[1mm]
  &\Upsilon_{mplk}=\mathcal{E}_{mi_2l}^{k}+\Delta y\sum_{q=1}^{\iota}\overline{a}_{pq}\delta_{y}\Upsilon_{mqlk},\quad \forall~i_2=0,1,\cdots,I_2\label{GME_sub5}\\[1mm]
  &\mathcal{E}_{m(i_{2}+1)l}^{k}=\mathcal{E}_{mi_2l}^{k}+\Delta y\sum_{p=1}^{\iota}\overline{b}_{p}\delta_{y}\Upsilon_{mplk},\quad \forall~i_2=0,1,\cdots,I_2\label{GME_sub6}\\[1mm]
  &\Upsilon_{mplk}=\mathcal{E}_{mpi_3}^{k}+\Delta z\sum_{v=1}^{\sigma}\widehat{a}_{lv}\delta_{z}\Upsilon_{mpv k},\quad \forall~i_3=0,1,\cdots,I_3\label{GME_sub7}\\[1mm]
  &\mathcal{E}_{mp(i_3+1)}^{k}=\mathcal{E}_{mpi_3}^{k}+\Delta z\sum_{u=1}^{\sigma}\widehat{b}_{u}\delta_{z}\Upsilon_{mpuk},\quad \forall~i_3=0,1,\cdots,I_3\label{GME_sub8}\\[1mm]
  &\tau{\rm {\bf K}}\delta_{t}\Upsilon_{mplk}+\tau\sum_{i=1}^3{\rm {\bf L}}_{i}\delta_{x_i}\Upsilon_{mplk}=\nabla S(\Upsilon_{mplk})\Delta W_m^k.\label{GME_sub9}
  \end{align}
\end{subequations}

\begin{thm}
Assume that the coefficients $a_{kj},\widetilde{a}_{mn},\overline{a}_{pq},\widehat{a}_{lv},b_k,\widetilde{b}_m,\overline{b}_{p},\widehat{b}_{u}$ of \eqref{RK00} satisfy the relations
\begin{align}\label{cond1}
\begin{split}
  &b_kb_j-b_ka_{kj}-b_ja_{jk}=0,\qquad \overline{b}_p\overline{b}_q-\overline{b}_p\overline{a}_{pq}-\overline{b}_q\overline a_{qp}=0\\
  &\widetilde{b}_m\widetilde{b}_n-\widetilde{b}_m\widetilde{a}_{mn}-\widetilde{b}_n\widetilde{a}_{nm}=0,\qquad \widehat{b}_l\widehat{b}_v-\widehat{b}_l\widehat{a}_{lv}-\widehat{b}_v\widehat{a}_{lv}=0
\end{split}
\end{align}
for all $k,j=1,\cdots,r$, $n,m=1,\cdots,s$, $p,q=1,\cdots,\iota$ and $l,v=1,\cdots,\sigma,$ then the stochastic RK methods \eqref{RK00} are stochastic multi-symplectic with the discrete stochastic multi-symplectic conservation law
\begin{equation}\label{def_sm_rk_maxwell}
  \frac{\omega^{\rho+1}-\omega^{\rho}}{\tau}+\frac{\kappa_{i_1+1}^{(1)}-\kappa_{i_1}^{(1)}}{\Delta x}+\frac{\kappa_{i_2+1}^{(2)}-\kappa_{i_2}^{(2)}}{\Delta y}+\frac{\kappa_{i_3+1}^{(3)}-\kappa_{i_3}^{(3)}}{\Delta z}=0,~a.s.,
\end{equation}
where
\begin{equation*}
\begin{split}
\omega^\rho&=\frac{1}{2}\sum_{m=1}^{s}\sum_{p=1}^{\iota}\sum_{l=1}^{\sigma}\widetilde{b}_m\overline{b}_p\widehat{b}_l
  {\rm d}z_{mpl}^{\rho}\wedge {\rm {\bf K}}{\rm d}z_{mpl}^{\rho},\quad \rho=0,1,\cdots,N,\\
\kappa_{i_1}^{(1)}&=\frac{1}{2}\sum_{k=1}^{r}\sum_{p=1}^{\iota}\sum_{l=1}^{\sigma}b_k\overline{b}_p\widehat{b}_l
  {\rm d}z_{i_1 pl}^{k}\wedge {\rm {\bf L}}_1{\rm d}z_{i_1 pl}^{k},\quad i_1=0,1,\cdots,I_1,\\
\kappa_{i_2}^{(2)}&=\frac{1}{2}\sum_{k=1}^{r}\sum_{m=1}^{s}\sum_{l=1}^{\sigma}b_k\widetilde{b}_m\widehat{b}_l
  {\rm d}z_{mi_2 l}^{k}\wedge {\rm {\bf L}}_2{\rm d}z_{mi_2 l}^{k},\quad i_2=0,1,\cdots,I_2,\\
\kappa_{i_3}^{(3)}&=\frac{1}{2}\sum_{k=1}^{r}\sum_{m=1}^{s}\sum_{p=1}^{\iota}b_k\widetilde{b}_m\overline{b}_p
  {\rm d}z_{mpi_3}^{k}\wedge {\rm {\bf L}}_3{\rm d}z_{mpi_3}^{k},\quad i_3=0,1,\cdots,I_3.
  \end{split}
\end{equation*}
\end{thm}
\begin{pr}
The proof is
analogous to that of the Theorem \ref{cond0}, so we ignore it here.
\end{pr}

For stochastic Maxwell equations \eqref{stochastic maxwell equations}, \cite{HJZC2017} presents that the energy is invariant, i.e.,
\begin{equation*}
  I(t):=\int_{D}(|\textbf{E}(x,y,z,t)|^{2}+|\textbf{H}(x,y,z,t)|^{2})dxdydz=Constant,~a.s.
\end{equation*}
In the context of our methods, it is interesting to see whether this energy functional $I(t)$ remain invariant, thus describing
the persistence of integrability of the fully discrete scheme. The result is stated in the following theorem.

\begin{thm}\label{Theorem}
Under the periodic boundary condition and the conditions \eqref{cond1}, the stochastic Runge-Kutta methods \eqref{RK00} have the following discrete energy conservation law, that is, for all $\rho=0,1,\cdots,N$
\begin{align}\label{58}
\begin{split}
\sum_{i_1}&\sum_{i_2}\sum_{i_3}\sum_{m=1}^{s}\sum_{p=1}^{\iota}\sum_{l=1}^{\sigma}\widetilde{b}_{m}\overline{b}_p\widehat{b}_l
\Big(|\mathbf{E}_{mpl;i_1,i_2,i_3}^{\rho+1}|^2+|\mathbf{H}_{mpl;i_1,i_2,i_3}^{\rho+1}|^2\Big)\\
&=\sum_{i_1}\sum_{i_2}\sum_{i_3}\sum_{m=1}^{s}\sum_{p=1}^{\iota}\sum_{l=1}^{\sigma}\widehat{b}_{m}\overline{b}_p\widehat{b}_l
\Big(|\mathbf{E}_{mpl;i_1,i_2,i_3}^{\rho}|^2+|\mathbf{H}_{mpl;i_1,i_2,i_3}^{\rho}|^2\Big),~a.s.
\end{split}
\end{align}
\end{thm}

\begin{pr}
Define energy functional $\zeta(\mathcal{E})=\mathcal{E}^T\mathcal{E}$, we have
 \begin{equation*}\label{Proof_6}
 \begin{split}
 \zeta(\mathcal{E}_{mpl}^{\rho+1})- \zeta(\mathcal{E}_{mpl}^{\rho})&=\left(\mathcal{E}_{mpl}^{\rho+1}\right)^T\mathcal{E}_{mpl}^{\rho+1}-\left(\mathcal{E}_{mpl}^{\rho}\right)^T\mathcal{E}_{mpl}^{\rho}\\
 &=\left(\mathcal{E}_{mpl}^{\rho+1}\right)^T\left(\mathcal{E}_{mpl}^{\rho+1}-\mathcal{E}_{mpl}^{\rho}\right)+\left(\mathcal{E}_{mpl}^{\rho+1}-\mathcal{E}_{mpl}^{\rho}\right)^T\mathcal{E}_{mpl}^{\rho}. \end{split}
 \end{equation*}
It follows from \eqref{GME_sub1} and \eqref{GME_sub2} that
 \begin{equation}\label{Proof_7}
 \begin{split}
 &\zeta(\mathcal{E}_{mpl}^{\rho+1})- \zeta(\mathcal{E}_{mpl}^{\rho})=\tau\left(\mathcal{E}_{mpl}^{\rho+1}\right)^T\sum_{k=1}^{r}b_k\delta_t\Upsilon_{mplk}+\tau\sum_{k=1}^{r}b_k\left(\delta_t\Upsilon_{mplk}\right)^T\mathcal{E}_{mpl}^{\rho}\\
 &=\tau\left(\mathcal{E}_{mpl}^{\rho}+\tau\sum_{k=1}^{r}b_k\delta_t\Upsilon_{mplk}\right)^T\sum_{k=1}^{r}b_k\delta_t\Upsilon_{mplk}+\tau\sum_{k=1}^{r}b_k\left(\delta_t\Upsilon_{mplk}\right)^T\mathcal{E}_{mpl}^{\rho}\\
&=\tau\sum_{k=1}^{r}b_k\left((\mathcal{E}_{mpl}^{\rho})^T\delta_t\Upsilon_{mplk}+(\delta_t\Upsilon_{mplk})^T\mathcal{E}_{mpl}^{\rho}\right)+\tau^2\sum_{k=1}^{r}\sum_{j=1}^{r}b_kb_j(\delta_t\Upsilon_{mplk})^T(\delta_t\Upsilon_{mplj})\\
&=\tau\sum_{k=1}^{r}b_k\left(\left(\Upsilon_{mplk}-\tau\sum_{j=1}^{r}a_{k_j}\delta_t \Upsilon_{mplj}\right)^T\delta_t\Upsilon_{mplk}+\left(\delta_t\Upsilon_{mplk}\right)^T\left(\Upsilon_{mplk}-\tau\sum_{j=1}^{r}a_{k_j}\delta_t \Upsilon_{mplj}\right)\right)\\
&~~~~~~~~~~~~+\tau^2\sum_{k=1}^{r}\sum_{j=1}^{r}b_kb_j(\delta_t\Upsilon_{mplk})^T(\delta_t\Upsilon_{mplj})\\
&=2\tau\sum_{k=1}^{r}b_k\Upsilon_{mplk}^T\delta_t\Upsilon_{mplk}+\tau^2\sum_{k=1}^{r}\sum_{j=1}^{r}\left(b_kb_j-b_ka_{kj}-b_ja_{jk}\right)(\delta_t\Upsilon_{mplk})^T(\delta_t\Upsilon_{mplj})\\
&=2\tau\sum_{k=1}^{r}b_k\Upsilon_{mplk}^T\delta_t\Upsilon_{mplk},
 \end{split}
 \end{equation}
where the last equality is due to \eqref{cond1}.

Denote $\Lambda_i:={\rm {\bf K}}^{-1}{\rm {\bf L}}_i,~i=1,2,3$ and define $\alpha_i(\mathcal{E})=\mathcal{E}^T\Lambda_i\mathcal{E}$, we can obtain
\begin{equation*}\label{Proof_8}
 \begin{split}
 \alpha_1(\mathcal{E}_{(i_1+1)pl}^{\rho})&- \alpha_1(\mathcal{E}_{i_1pl}^{\rho})=\left(\mathcal{E}_{(i_1+1)pl}^{\rho}\right)^T\Lambda_1\mathcal{E}_{(i_1+1)pl}^{\rho}-\left(\mathcal{E}_{i_1pl}^{\rho}\right)^T\Lambda_1\mathcal{E}_{i_1pl}^{\rho}\\
 &=\left(\mathcal{E}_{(i_1+1)pl}^{\rho}\right)^T\Lambda_1\left(\mathcal{E}_{(i_1+1)pl}^{\rho}-\mathcal{E}_{i_1pl}^{\rho}\right)+\left(\mathcal{E}_{(i_1+1)pl}^{\rho}-\mathcal{E}_{i_1pl}^{\rho}\right)^T\Lambda_1\mathcal{E}_{i_1pl}^{\rho}. \end{split}
 \end{equation*}
It follows from \eqref{GME_sub3} and \eqref{GME_sub4} that
 \begin{equation}\label{Proof_9}
 \begin{split}
 &\alpha_1(\mathcal{E}_{(i_1+1)pl}^{\rho})- \alpha_1(\mathcal{E}_{i_1pl}^{\rho})\\
 &=\Delta x\left(\mathcal{E}_{(i_1+1)pl}^{\rho}\right)^T\Lambda_1\sum_{m=1}^{s}\tilde{b}_m\delta_x\Upsilon_{mplk}+\Delta x\sum_{m=1}^{s}\tilde{b}_m\left(\delta_x\Upsilon_{mplk}\right)^T\Lambda_1\mathcal{E}_{i_1pl}^{\rho}\\
 &=\Delta x\left(\mathcal{E}_{i_1pl}^{\rho}+\Delta x\sum_{m=1}^{s}\tilde{b}_m\delta_x\Upsilon_{mplk}\right)^T\Lambda_1\sum_{m=1}^{s}\tilde{b}_m\delta_x\Upsilon_{mplk}+\Delta x\sum_{m=1}^{s}\tilde{b}_m\left(\delta_x\Upsilon_{mplk}\right)^T\Lambda_1\mathcal{E}_{i_1pl}^{\rho}\\
&=\Delta x\sum_{m=1}^{s}\tilde{b}_m\left((\mathcal{E}_{i_1pl}^{\rho})^T\Lambda_1\delta_x\Upsilon_{mplk}+(\delta_x\Upsilon_{mplk})^T\Lambda_1\mathcal{E}_{i_1pl}^{\rho}\right)+\Delta x^2\sum_{m=1}^{s}\sum_{n=1}^{s}\tilde{b}_m\tilde{b}_n(\delta_x\Upsilon_{mplk})^T\Lambda_1\delta_x\Upsilon_{nplk}\\
&=\Delta x\sum_{m=1}^{s}\tilde{b}_m\Bigg(\left(\Upsilon_{mplk}-\Delta x\sum_{n=1}^{s}a_{k_j}\delta_x \Upsilon_{nplk}\right)^T\Lambda_1\delta_x\Upsilon_{mplk}\\
&+\left(\delta_x\Upsilon_{mplk}\right)^T\Lambda_1\left(\Upsilon_{mplk}-\Delta x\sum_{n=1}^{s}a_{k_j}\delta_x \Upsilon_{nplk}\right)\Bigg)+\Delta x^2\sum_{m=1}^{s}\sum_{n=1}^{s}\tilde{b}_mb_j(\delta_x\Upsilon_{mplk})^T\Lambda_1(\delta_x\Upsilon_{nplk})\\
&=2\Delta x\sum_{m=1}^{s}\tilde{b}_m\Upsilon_{mplk}^T\Lambda_1\delta_x\Upsilon_{mplk}+\Delta x^2\sum_{m=1}^{s}\sum_{n=1}^{s}\left(\tilde{b}_m\tilde{b}_n-\tilde{b}_m\tilde{a}_{mn}-\tilde{b}_n\tilde{a}_{nm}\right)(\delta_x\Upsilon_{mplk})^T\Lambda_1(\delta_x\Upsilon_{nplk})\\
&=2\Delta x\sum_{m=1}^{s}\tilde{b}_m\Upsilon_{mplk}^T\Lambda_1\delta_x\Upsilon_{mplk},
 \end{split}
 \end{equation}
where the last two equalities are due to the symmetry of matrix $\Lambda_1$ and the multi-symplectic conditions \eqref{cond1}, respectively.

Similarly, we can prove that
\begin{equation}\label{Proof_10}
 \begin{split}
 \alpha_2(\mathcal{E}_{m(i_2+1)l}^{\rho})- \alpha_2(\mathcal{E}_{mi_2l}^{\rho})=2\Delta y\sum_{p=1}^{\iota}\bar{b}\Upsilon_{mplk}^T\Lambda_2\delta_y\Upsilon_{mplk},
 \end{split}
 \end{equation}
\begin{equation}\label{Proof_11}
 \begin{split}
 \alpha_3(\mathcal{E}_{mp(i_3+1)}^{\rho})- \alpha_3(\mathcal{E}_{mpi_3}^{\rho})=2\Delta z\sum_{l=1}^{\sigma}\hat{b}\Upsilon_{mplk}^T\Lambda_3\delta_z\Upsilon_{mplk}.
 \end{split}
 \end{equation}

Next, by using the non-singularity of matrix ${\rm {\bf K}}$ and the equality \eqref{SS}, multiplying both sides of \eqref{GME_sub9} with $\Upsilon_{mplk}^T{\rm {\bf K}}^{-1}$, it yields
\begin{align}\label{variable}
\begin{split}
  \Upsilon_{mplk}^T\delta_t\Upsilon_{mplk}+\Upsilon_{mplk}^T\sum_{i=1}^3\Lambda_i\delta_{x_i}\Upsilon_{mplk}=\lambda\Upsilon_{mplk}^T{\rm {\bf K}}^{-1}\Upsilon_{mplk}\frac{\Delta W_m^k}{\tau}.
\end{split}
\end{align}
Due to the skew-symmetry of matrix ${\rm {\bf K}}$, summing up for $m,p,l,k$, we have
\begin{equation*}\label{111}
\sum_{k=1}^{r}\sum_{m=1}^{s}\sum_{p=1}^{\iota}\sum_{l=1}^{\sigma}b_k\tilde{b}_m\bar{b}_p\hat{b}_l\left(\Upsilon_{mplk}^T\delta_t\Upsilon_{mplk}+\Upsilon_{mplk}^T\sum_{i=1}^3\Lambda_i\delta_{x_i}\Upsilon_{mplk}\right)=0,
\end{equation*}
then, multiplying the above equation by $2\tau\Delta x\Delta y\Delta z$ and substituting \eqref{Proof_7}, \eqref{Proof_9}, \eqref{Proof_10} and \eqref{Proof_11} into the above equation, it holds
\begin{align*}
\begin{split}
\Bigg(\Delta x\Delta y\Delta z&\sum_{m=1}^{s}\sum_{p=1}^{\iota}\sum_{l=1}^{\sigma}\tilde{b}_m\bar{b}_p\hat{b}_l\left(\zeta(\mathcal{E}_{mpl}^{\rho+1})- \zeta(\mathcal{E}_{mpl}^{\rho})\right)\\
&+\tau\Delta y\Delta z\sum_{k=1}^{r}\sum_{p=1}^{\iota}\sum_{l=1}^{\sigma}b_k\bar{b}_p\hat{b}_l\left(\alpha_1(\mathcal{E}_{(i_1+1)pl}^{\rho})- \alpha_1(\mathcal{E}_{i_1pl}^{\rho})\right)\\
&+\tau\Delta x\Delta z\sum_{k=1}^{r}\sum_{m=1}^{s}\sum_{l=1}^{\sigma}b_k\tilde{b}_m\hat{b}_l\left(\alpha_2(\mathcal{E}_{m(i_2+1)l}^{\rho})- \alpha_2(\mathcal{E}_{mi_2l}^{\rho})\right)\\
&+\tau\Delta x\Delta y\sum_{k=1}^{r}\sum_{m=1}^{s}\sum_{p=1}^{\iota}b_k\tilde{b}_m\bar{b}_p\left(\alpha_3(\mathcal{E}_{mp(i_3+1)}^{\rho})- \alpha_3(\mathcal{E}_{mpi_3}^{\rho})\right)
\Bigg)=0.
\end{split}
\end{align*}
Summing up for $i_1, i_2, i_3$ over the spatial domain and utilizing the periodic boundary condition, we can get the result \eqref{58}.
 \end{pr}

\begin{rmk}
The results of these two theorems are evidently consistent with the stochastic multi-symplectic conservation law and discrete energy conservation law in \cite{HJZC2017}, respectively, which means that the stochastic multi-symplectic conservation law and energy conservation law can be exactly preserved by the proposed stochastic multi-symplectic RK methods \eqref{RK00}-\eqref{cond1}.
\end{rmk}

For the general stochastic Hamiltonian PDEs in the sense of the Stratonovich
\begin{align}
{\rm {\bf K}}{\rm d}_tz+\sum_{m=1}^{M}{\rm {\bf L}}_mz_{x_m}{\rm d}t=\nabla_z S_1(z){\rm d}t+\nabla_z S_2(z)\circ{\rm d}W(t),~~z\in\mathbb{R}^n,
\end{align}
where ${\rm {\bf K}}$ and ${\rm {\bf L}}_m,~m=1,2,\cdots,M$ are skew-symmetric matrices, we can obtain the following corollary.
\begin{corr}
Let the matrix ${\rm {\bf K}}$ be nonsingular and the Hamiltonian functions $S_1=\frac{1}{2}z^T{\rm {\bf A}}z,~S_2=\frac{1}{2}z^T{\rm {\bf B}}z$, where ${\rm {\bf A}},~{\rm {\bf B}}$ are symmetric matrices. If ${\rm {\bf K}}^{-1}{\rm {\bf A}}$ and ${\rm {\bf K}}^{-1}{\rm {\bf B}}$ are skew-symmetric matrices, then for all $\rho=0,1,\cdots,N$, it holds
\begin{align}
\begin{split}
\sum_{i_1}\cdots\sum_{i_M}&\sum_{j_1=1}^{s_1}\cdots\sum_{j_M=1}^{s_M}b_{j_1}^{(1)}\cdots b_{j_M}^{(M)}|z_{j_1\cdots j_M;i_1,\cdots,i_M}^{\rho+1}|^2\\
&=\sum_{i_1}\cdots\sum_{i_M}\sum_{j_1=1}^{s_1}\cdots\sum_{j_M=1}^{s_M}b_{j_1}^{(1)}\cdots b_{j_M}^{(M)}|z_{j_1\cdots j_M;i_1,\cdots,i_M}^{\rho}|^2,~a.s.,
\end{split}
\end{align}
where $i_1,\cdots,i_M$ denote the spatial grids, $s_1,\cdots,s_M$ denote the stage of stochastic RK methods applied to the directions of space, $|z|^2$ denotes the sum of components, i.e., $|z|^2:=|z_1|^2+\cdots+|z_d|^2$.
\end{corr}

\section{Conclusions}
In this paper, we construct a class of stochastic RK methods to stochastic Hamiltonian PDEs and give some sufficient conditions for stochastic multi-symplecticity of the constructed stochastic RK methods. Finally, we apply these techniques to three-dimensional stochastic Maxwell equations with multiplicative noise in sense of Stratonovich. We show that the proposed stochastic RK methods can preserve the discrete stochastic multi-symplectic conservation law and the discrete energy conservation law almost surely. However, the theoretical analysis of the convergence of stochastic multi-symplectic RK methods is difficult, we will devote to study it rigorously in the future work.






\bibliography{ref.bib}
\section{References}







\end{document}